\documentclass{article}

\usepackage{amsfonts}
\usepackage{color}

\definecolor{Red}{cmyk}{0,1,1,0}

\definecolor{Blue}{cmyk}{1,1,0,0}

\definecolor{Pink}{cmyk}{0,1,0,0}

\definecolor{Green}{cmyk}{1,0,1,0.5}

\newcommand{\ba}{\begin{array}}
\newcommand{\ea}{\end{array}}
\newcommand{\be}{\begin{equation}}
\newcommand{\ee}{\end{equation}}
\newcommand{\ben}{\begin{enumerate}}
\newcommand{\een}{\end{enumerate}}

\newcommand{\eop}{\hfill \rule{0.7ex}{1.6ex}}

\newcommand{\F}{ {\cal F} }

\newcommand{\R}{\mathbb{R}}

\newcommand{\B}{{\cal B}}

\newcommand{\N}{\mathbb{N}}

\newcommand{\dlim}{\displaystyle \lim}

\newcommand{\dsum}{\displaystyle \sum}
\newcommand{\dsup}{\displaystyle \sup}

\newcommand{\dmin}{\displaystyle \min}

\newcommand{\oo}{\infty}

\newtheorem{teo}{Theorem}
\newtheorem{lema}{Lemma}
\newtheorem{cor}{Corollary}

\title{
{\large{ \bf{RENORMALIZATION GROUP ANALYSIS OF NONLINEAR DIFFUSION
EQUATIONS WITH TIME DEPENDENT COEFFICIENTS: ANALYTICAL RESULTS}}}
\footnotetext{
\indent \,\,\ This work is partially supported by CNPq-Brazil.} }

\author{
{Gast\~ao A. Braga$^1$, Frederico Furtado$^2$},
\\
{Jussara M. Moreira$^1$},
{Leonardo T. Rolla$^3$}\\
\\
$^1$\small{Departamento de Matem\'atica -
Universidade Federal de Minas Gerais}\\
\small{Caixa Postal 1621,} \small{30161-970 -
Belo Horizonte - MG - Brazil}\\\\
$^2$\small{Department of Mathematics - University of Wyoming}\\
\small{Laramie - 82071 - USA}\\\\
$^3$\small{Instituto de Matem\'atica Pura e Aplicada}\\
\small{Estrada Dona Castorina 110,} \small{22460-320 - Rio de}
Janeiro - RJ - Brazil}
\date{}
\begin{document}

\baselineskip = 14pt

\maketitle

\begin{abstract}
We study the long-time asymptotics of a certain class of nonlinear
diffusion equations with time-dependent diffusion coefficients which
arise, for instance, in the study of transport by randomly
fluctuating velocity fields. Our primary goal is to understand the
interplay between anomalous diffusion and nonlinearity in
determining the long-time behavior of solutions. The analysis
employs the renormalization group method to establish the
self-similarity and to uncover universality in the way solutions
decay to zero.
\end{abstract}

This preprint has the same numbering of sections, equations and theorems as the published article ``\emph{Discrete Contin. Dyn. Syst. Ser. B 7 (2007), 699--715}.''



\section{\large{Introduction}}
\label{sec:intr}

Theories of transport by a random velocity field are used in a
number of important problems in many fields of science and
engineering. Examples range from mass and heat transport in
geophysical flows \cite{bib:davis91}, to combustion and
chemical engineering \cite{bib:moffatt83}, to hydrology
\cite{bib:dagan87} and petroleum engineering \cite{bib:lake86}. In
each of these examples one finds important physical processes which
involve the transport of a passive scalar quantity in the presence
of a complex velocity field that fluctuates randomly on length and
time scales comparable to those on which the transport process
occurs. A central goal of transport theories is to understand the
effects produced by such random velocity fluctuations on the mean
(ensemble averaged) transport.

A characteristic feature of these theories, which are usually based
on perturbation methods, are infrared (long wavelength or low
frequency) divergences in the terms of the pertinent perturbation
expansions due to long-range (nonintegrable) spatial or temporal
correlations in velocity fluctuations. See, e.g., \cite{bib:avell,
bib:glimm, bib:majda} and references therein. Typically, the
dominant divergences in lowest order are of diffusion type, and
correspond to anomalous diffusion in the ensemble averaged transport
equations, for which the effective diffusion coefficient increases
with time and is divergent as time $t \to \infty$. So, according to
these theories the mean concentration, $u$, of a passive scalar
field being advected by a random velocity field with strong,
long-range correlations satisfies, under appropriate conditions, an
equation of the form \cite{bib:glimm, bib:majda} \be \label{eq:form}
u_t=c(t)u_{xx} + F(u, u_x, u_{xx}), \quad c(t) \sim t^p \mbox{ as }
t\to\infty, \mbox{ with } p > 0. \ee It is our intent here and in
\cite{bib:gas-fred-ju-leo-tp1} to analyze the long-time behavior of
solutions of equations of the form~(\ref{eq:form}). The analysis
assumes $F = F(u)$ to be superlinear, in the sense that $F(u) =
O(u^\alpha)$ as $u\to 0$, with $\alpha >1$. The inclusion of the
nonlinear term $F(u)$ in~(\ref{eq:form}) accounts for situations in
which the scalar field is not conservative, meaning that its
concentration $u$ undergoes changes due to physical, chemical or
biological processes.

We are concerned primarily with the interplay between anomalous
diffusion (measured in terms of the exponent $p$) and nonlinearity
(measured in terms of the exponent $\alpha$) in determining the
scaling behavior of solutions as they decay to zero. In the present
paper we analyze the situation where the nonlinearity is analytic
and ``supercritical'' (or irrelevant), in the sense that $F(u) =
\sum_{j\geq \alpha} a_j u^j$, with $\alpha > (p+3)/(p+1)$.  We show
that in this case diffusion is the dominant effect in the limit
$t\to\infty$, and determines the scaling form of solutions with
sufficiently localized initial data as they decay to zero: \be
\label{eq:asym-limi} u(x,t) \sim t^{-\gamma/2} \phi\left(
\frac{x}{t^{(p+1)/2}}\right) \quad \mbox{as } t\to\infty, \ee where
$\gamma = p+1$ and the function on the right-hand side
of~(\ref{eq:asym-limi}) is a self-similar solution of
equation~(\ref{eq:form}) with $c(t) = t^p$ and $F\equiv 0$. Thus, a
curious phenomenon, \emph{universality}, is characterized: solutions
of many different equations, but differing only in the nonlinear
term $F(u)$, the higher-order asymptotics of $c(t)$, or both,
nevertheless share the same asymptotic scaling behavior, given by a
self-similar solution of the time-dependent diffusion equation $u_t
= t^p\, u_{xx}$.  Thus, it can be said that such equations belong to
the same \emph{universality class}, in that all of the members of
this class exhibit the same asymptotic behavior, insofar as the
scaling behavior of their solutions (with initial data in suitable
classes) is the same.

``Subcritical'' (or relevant) and ``critical'' (or marginal)
nonlinearities, namely those varying as $F(u)\sim u^\alpha$ with
$\alpha < (p+3)/(p+1)$ and $\alpha = (p+3)/(p+1)$, respectively, in
the limit $u\to 0$, are analyzed numerically in
\cite{bib:gas-fred-ju-leo-tp1}. This analysis does not assume $F(u)$
to be analytic at $u = 0$. In marked contrast with the supercritical
case, in the subcritical case the asymptotic scaling behavior of
solutions is strongly affected by the nonlinear term $F(u)$.  In
particular, the decay exponent $\gamma = 2/(\alpha - 1)$ is
determined by the leading-order term ($u^\alpha$) in $F(u)$, and the
function on the right-hand side of~(\ref{eq:asym-limi}) is now a
self-similar solution of the time-dependent reaction-diffusion
equation $u_t = t^p\, u_{xx} - u^\alpha$. Thus, the phenomenon of
universality is again visible. Equations differing only in the
higher-order asymptotics of $F(u)$, $c(t)$, or both, fall in the
same universality class.

The critical case is peculiar.
It marks the crossover from a scaling regime controlled (mostly) by
diffusion (supercritical case) to a scaling regime strongly
influenced, and in certain aspects determined, by nonlinearity
(subcritical case). Thus, in the critical case neither diffusion nor
nonlinearity prevails, and the scaling regime which is observed
bears some features of the supercritical one (same $\gamma$ and
scaling function $\phi$), but acquires an extra logarithmic decay
factor, the imprint of the critical nonlinearity.

The following heuristic arguments motivate the critical case, $\alpha=(p+3)/(p+1)$.
Anticipating that solutions $u\to 0$ as $t\to\infty$, it is possible to
simplify the analysis by taking $F(u)=\lambda u^\alpha$. This can be thought of as
being an approximation near $u=0$. With this choice, and recognizing that it is the
large-time regime which we seek to understand, we simplify the equation to
\[
u_t = t^p u_{xx} + \lambda u^\alpha.
\]
Under the parabolic scaling $x\mapsto L^{(p+1)/2} x$, $t\mapsto Lt$,
$u\mapsto L^{(p+1)/2} u$, $L\gg1$, this equation becomes
\[
u_t = t^p u_{xx} + L^{1+(1-\alpha)(p+1)/2}\lambda u^\alpha.
\]
Thus, provided $\alpha > (p+3)/(p+1)$, as we iterate such scaling transformation
we end up with an equation where the nonlinear term decreases by a factor
$L^{1+(1-\alpha)(p+1)/2}$ at each rescaling. Consequently, as the number of rescalings
$n\to\infty$ (equivalently, $t\to\infty$), the linear diffusion term dominates the nonlinear
term and we may expect solutions to decay as $t^{-(p+1)/2}$, the rate determined by the
linear diffusion. This argument fails when $\alpha=(p+3)/(p+1)$, which we dub the critical
case.

We now state precisely the main result of this paper. For this
purpose we introduce the spaces
$$
{\cal B}_q \equiv \{f:\R \rightarrow \R \,\ | \,\ \hat{f}(w) \in
C^1(\R) \ {\mbox{and}} \ \|f\|< \oo \}, \quad q > 1,
$$
with norm $\|f\| = \sup_{w \in \R}\left[(1+|w|^q)\left(\left|\hat{f}(w)\right|+
\left|\hat{f}'(w)\right|\right)\right]$ and
$$
B^{(\infty)}\equiv \{u:\R \times [1,+\infty) \rightarrow \R \,\ |
\,\ u(\cdot,t)\in \B_q \mbox{ for all }  t\geq1 \ {\mbox{and}} \
\|u\|_{\infty}< \oo \},
$$
where  $\|u\|_{\infty}=\sup_{t \geq 1}\|u(\cdot,t)\|$.

Consider the following initial value problem (IVP)
\be
\label{pvi:eqc(t)F}
\cases{u_t=c(t)u_{xx}+\lambda F(u), & $t>1,
\,\,\ x\in\R$ \cr \cr u(x,1)=f(x), & $x\in\R$}
\ee and assumptions:
\begin{itemize}
\item[($H1$)] \label{hyp:h1} $f\in {\cal B}_q$ for some $q>1$;
\item[($H2$)] \label{hyp:h2}  $c(t)>0$ for $t>1$ and $c(t) = t^p + o(t^p)$
as $t\to\infty$, with $p>0$;
\item[($H3$)] \label{hyp:h3} $F(u) = \sum_{j\geq \alpha} a_j u^j$
analytic at $u=0$, with $\alpha > (p+3)/(p+1)$.
\end{itemize}
We shall prove the following.
\begin{teo}
\label{teo:compassgeral} Assume $(H1)-(H3)$. Then there exists an
$\varepsilon>0$ such that, for  $\|f\|<\varepsilon$, one can find
$B \subset B^{(\infty)}$ such that the IVP~(\ref{pvi:eqc(t)F}) has a unique
solution $u \in B$ which satisfies, for some constant $A$,
\be
\label{eq:cotanaolin} \dlim_{t \to \oo}
\left\|\sqrt{t^{p+1}}u(\sqrt{t^{p+1}}\cdot,t)-Af_p^*(\cdot)\right\|=0,
\ee with
$f_p^*(x)=\sqrt{\frac{p+1}{4\pi}}\,\,e^{-\frac{(p+1)}{4}x^2}$.
\end{teo}

Our proof relies on the Renormalization Group (RG) approach. RG
methods were originally introduced, and proved to be very useful, in
quantum field theory \cite{bib:bogoliubov,bib:gell-low} and
statistical mechanics \cite{bib:wilson71a,bib:wilson71b}. Their application to the
asymptotic analysis of deterministic differential equations (both
ODEs and PDEs) was initiated and developed by Goldenfeld, Oono and
collaborators \cite{bib:gold3, bib:gold2, bib:gold1}.  See
\cite{bib:gold-book} and \cite{bib:oono} for detailed accounts. The
mathematical aspects of the method were rigorously established by
Bricmont, Kupiainen and collaborators \cite{bib:bric-kupa,
bib:bric-kupa-lin}. See also \cite{bib:cag, bib:merd}.

In the RG approach the long-time behavior of solutions to PDEs is
related to the existence and stability of fixed points of an
appropriate RG transformation. The definition of an RG
transformation involves two basic steps.  The first step is the
integration (solution) of the PDE over a finite time interval; its
purpose is to eliminate the ``small time'' information in the
problem (coarse graining). The second step is rescaling, to change
the time scales in proportion to those eliminated (by integration),
so that a nominally constant time scale is under study.  The
iterative application of the RG transformation progressively evolves
the solution in time and at the same time renormalizes the terms of
the PDE. Once a
proper RG transformation has been found for a particular problem,
these terms are divided into two types: neutral and irrelevant,
according to whether their magnitude is unchanged or decreases with
each RG iteration. The irrelevant terms iterate to zero and the
dynamics at large times is then controlled by the neutral terms.
This accounts for the observed universal scaling behavior of
solutions as they decay to zero. Thus, the RG provides a natural
framework in which to understand universality.

Our results contribute to a large body of literature devoted to the
study of long-time asymptotics of nonlinear PDEs.  The work of
Barenblatt and collaborators \cite{bib:barenblatt-book} has had a
major impact in this field of study, specially in elucidating the
importance, as well as intricacies, of self-similarity in
intermediate asymptotics. Our analysis follows closely the one in
\cite{bib:bric-kupa-lin}. See also \cite{bib:bona1,
bib:bric-kup-xin, bib:li-qi}. In spirit, our
contribution relates also to the work summarized in
\cite{bib:aronson-vazquez,bib:vazquez-ima}.

Our analysis can easily be extended in a number of directions.
These include equations in more than one space dimension, nonlinearities
involving $u$ as well as derivatives of $u$ (as in
equation~(\ref{eq:form})), and nonlinearities with time-dependent coefficients.
The modifications needed in each case are straightforward. For instance, if the
nonlinearity is of the form $d(t)F(u)$, with $d(t)\sim t^r$ as $t\to\oo$ and
$F(u) \sim u^\alpha$ as $u\to 0$, the elementary scaling argument presented
above suggests that the critical exponent is now
$\alpha = (p+3+2r)/(p+1)$. With this proviso Theorem~\ref{teo:compassgeral}
should also hold in this case.

The rest of this paper is organized as follows.  In
Section~\ref{sec:exis-uniq} we establish the existence and
uniqueness of solutions for problem~(\ref{pvi:eqc(t)F}) over a
finite time interval. In Section~\ref{sec:asint} we employ the RG
approach to extend this result to an infinite time interval and
obtain the long-time asymptotics of solutions, thereby proving
Theorem~\ref{teo:compassgeral}.


\section{\large{Local Existence and Uniqueness}}
\label{sec:exis-uniq}

In this section we prove the existence and uniqueness of solutions
for problem (\ref{pvi:eqc(t)F}) over a finite time interval using a
fixed-point argument.  In the next section we employ the RG
iterative procedure \cite{bib:bric-kupa-lin} to extend this local
result over an infinite time interval. In the process we obtain
upper bounds that lead to the limit (\ref{eq:cotanaolin}).

As a preliminary remark, we note that generically the constants obtained
in this and in the next sections depend on $q>1$, the function $c(t)$  and
the coefficients $a_j$ of $F(\cdot)$. However, for simplicity of notation
we omit this dependence. Also, without loss of generality we assume
$\lambda\in[-1,1]$ so that the estimates obtained will be valid
uniformly with respect to $\lambda$.

We start with the definition of certain spaces and operators that
will be used throughout this paper.  For $L>1$, we
introduce the space
\be
\label{def:spacebl}
B^{(L)}\equiv \{u:\R
\times [1,L] \rightarrow \R \,\ | \,\ u(\cdot,t)\in \B_q \mbox{ for
all }t \in [1,L]\},
\ee
with norm $\|u\|_L=\dsup_{t \in
[1,L]}\|u(\cdot,t)\|$. Next, let
$$
  u_f(x,t)=\frac{1}{\sqrt{4\pi s(t)}} \int_\R
  e^{-\frac{(x-y)^2}{4s(t)}}f(y)\; \mathrm{d}y,
$$
\be
\label{def:N(u)}
N(u)(x,t)=\lambda \int_1^{t}\frac{1}{\sqrt{4\pi
[s(t)-s(\tau)]}} \int_\R
e^{-\frac{(x-y)^2}{4[s(t)-s(\tau)]}}F(u(y,\tau))\;\mathrm{d}y\mathrm{d}\tau,
\ee
with
$s(t)=\int_1^t c(v)\;\mathrm{d}v$, and define the operator
\be
\label{def:T(u)}
T(u)\equiv u_f+ N(u).
\ee

We shall prove that the operator $T$ has a unique fixed point.
This is equivalent to proving the existence and uniqueness of
solutions to IVP~(\ref{pvi:eqc(t)F}).
Specifically, we prove existence and uniqueness of solutions
in the ball $B_f$ defined below, provided $\|f\|$ is sufficiently small.
The introduction of $B_f$ is a necessary step since we
only stipulate the behavior of $F(u)$ near $u=0$.
So we define
\be
\label{spacebf}
B_f \equiv \left\{u\in B^{(L)}:\|u-u_f\|_L\leq \|f\|\right\}.
\ee

\begin{teo}
\label{teo:existunic} Assume $(H1)$-$(H3)$ and let $L>1$. Then there
is an $\varepsilon>0$ such that the IVP~(\ref{pvi:eqc(t)F}) has a
unique solution $u$ in $B_f$, if $\|f\|<\varepsilon$.
\end{teo}

The result follows immediately from Lemmas~\ref{lema:Nungeral}
and~\ref{lema:Nuvngeral}, since the former establishes that $T$ maps
$B_f$ into itself and the latter that $T$ is a contraction.
We note that Theorem~\ref{teo:existunic} also holds under the weaker
assumption $\alpha \geq 2$, instead of $\alpha>(p+3)/(p+1)$ as
stated in $(H3)$. However, the latter is necessary for the proof of
Theorem~\ref{teo:compassgeral} in the next section.

\begin{lema}
\label{lema:Nungeral} Let $L > 1$.
There is an $\varepsilon'>0$ such that
$\|N(u)\|_L<\|f\|$ for all $u\in B_f$, if $f\in\B_q$ and $\|f\|<
\varepsilon'$.
\end{lema}
\begin{lema}
\label{lema:Nuvngeral}
For each $L>1$
there is an $\varepsilon'' >0$ such that
$\|N(u)-N(v)\|_L< \frac{1}{2}\|u-v\|_L$ for all $u, v\in B_f$, if
$f\in\B_q$ and $\|f\|<\varepsilon''$.
\end{lema}

Before proving Lemmas~\ref{lema:Nungeral} and~\ref{lema:Nuvngeral},
we make some remarks and obtain estimates which are needed in the
proof and also in Section~\ref{sec:asint}.

Let us assume that the Taylor expansion
of $F(u)$ at $u = 0$ has a finite radius of convergence $\rho$ (the
case $\rho = \infty$ is less restrictive).
We argue that for $N(u)$ and $T(u)$ to be well defined it suffices
to require that $u \in B_f$ with $f\in{\cal B}_q$ such that
$\|f\|<[2C_q(1+\sqrt{s(L)})]^{-1}\rho$.
Indeed, we have to check that $|u(x,t)|<\rho$ for
all $x\in \R$ and $t \in [1,L]$ and this will follow by comparison of
different norms.
First, notice that $\B_q \subset L^1(\R) \cap L^2(\R) \cap L^{\oo}(\R)$
by elementary Fourier transform calculations. For instance, for all $x\in\R$
$$
|h(x)| \leq
\sup_x \frac{1}{2\pi} \int_{-\oo}^{\oo}|\widehat{h}(w)e^{iwx}|\;\mathrm{d}w
\leq \frac{1}{2\pi}\int_{\R}\frac{\|h\|}{1+|w|^q}\;\mathrm{d}w=
C_q\|h\|.
$$
So, if $u \in B^{(L)}$, we have $u(\cdot,t)\in \B_q$ and
$$
  \sup_x |u(x,t)| \leq C_q \|u\|_L
$$
for all $t\in [1,L]$. Now for $f \in \B_q$ and $u \in B_f$
we have
\be
\label{ineq:u-f}
\|u\|_L \leq 2 \left(1+\sqrt{s(L)}\right)\|f\|.
\ee

Define
$\varphi(\tau)\equiv s(t)-s(t-\tau)$ (recall $s(t)=\int_1^t c(v)\;\mathrm{d}v$)
and $J(w,t)\equiv \smallint_0^{t-1}
w\varphi(\tau)e^{-\varphi(\tau)w^2}\;\mathrm{d}\tau$ for $t\geq 1$
and $w\in\R$.
The following estimate holds:
\be
\label{cotaJ}
|J(w,t)|\leq (t-1)\sqrt{s(t)}.
\ee
Indeed,
taking $x=\sqrt{\varphi(\tau)}w$ and using that $xe^{-x^2}\leq 1$
for all $x \in \R$,
\[
|J(w,t)|=\left|\int_0^{t-1}\sqrt{\varphi(\tau)}xe^{-x^2}\;\mathrm{d}\tau
\right| \leq \int_0^{t-1}\sqrt{\varphi(\tau)}\;\mathrm{d}\tau.
\]
Since the integrand is a continuous function in $[0,t-1]$ and
$\varphi(\tau)$ is an increasing function (notice that $c(t)\geq
0$), we can bound the last integral by
$(t-1)\sqrt{\varphi(t-1)}=(t-1)\sqrt{s(t)-s(1)}$ and, since $s(1)=0$,
we have~(\ref{cotaJ}).

Also, if $q>1$ and $w\in \R$, \be \label{cotaI}
I(w)\equiv\int_{\R}\frac{1}{1+|x|^q}\cdot \frac{1}{1+|x-w|^q}\;\mathrm{d}x \leq
\frac{C}{1+|w|^q}, \ee where \be \label{def:C} C = C(q) =
(2^{q+1}+3)\int_{\R}\frac{1}{1+|x|^q}\;\mathrm{d}x. \ee Now, motivated by
hypothesis $(H3)$, let $C$ be given by~(\ref{def:C}) and define the
sums
\begin{eqnarray}
\label{somaS0} S_0(z) &\equiv& \sum_{j
\geq
\alpha}\left(\frac{C}{2\pi}\right)^{j-1}|a_j|z^j,\\
\label{somaS1} S_1(z)&\equiv& \sum_{j \geq
\alpha}\left(\frac{C}{2\pi}\right)^{j-1}|a_{j}|z^{j-2},\\
\label{somaS2} S_2(z)&\equiv& \sum_{j\geq
\alpha}\left(\frac{C}{2\pi}\right)^{j-1}j|a_j|z^{j-2}.
\end{eqnarray}
Notice that the radius of convergence of these sums is
$(2\pi\rho)/C<\rho$. We will now consider only those functions $u$
such that $|u(x,t)|<\pi\rho/C\equiv \rho_0$ for all $x\in \R$ and
$t \in [1,L]$. Invoking (\ref{ineq:u-f}), it thus suffices to take
$f$ such that
\be
\label{eq:normaf}
\|f\|<\left[2C_q(1+\sqrt{s(L)})\right]^{-1}\rho_0.
\ee


{\bf{Proof of Lemma~\ref{lema:Nungeral}:}}
Taking the Fourier
transform of $N(u)$ yields 
$$
\widehat{N(u)}(w,t) = \lambda \dsum_{j \geq
\alpha}a_j\int_0^{t-1}e^{-\varphi(\tau)
w^2}\widehat{u^j}(w,t-\tau)\;\mathrm{d}\tau.
$$
Writing $\widehat{u^j}$ as convolutions of $\hat u$, each term in
the sum above is of the form
\be
\label{eq:unchapeugeral}
\frac{a_j}{(2\pi)^{j-1}}\int_0^{t-1}\;\mathrm{d}\tau e^{-\varphi(\tau) w^2}
\int_{\R^{j-1}}\hat u(w-p_1)\hat u(p_1-p_2) \cdots \hat
u(p_{j-1})\;\mathrm{d}p_1\cdots \mathrm{d}p_{j-1},
\ee
where we have omitted the
dependence of $\hat u$ on $t-\tau$. Since the absolute value of
$\hat u$ is bounded by $\|u\|_L/(1+|w|^q)$, (\ref{eq:unchapeugeral})
can be upper-bounded by
$$
\frac{|a_j|}{(2\pi)^{j-1}}\|u\|_L^j \int_0^{t-1}\;\mathrm{d}\tau
e^{-\varphi(\tau) w^2}\int_{\R^{j-1}} \frac{1}{1+|w-p_1|^q} \cdots
\frac{1}{1+|p_{j-1}|^q}\;\mathrm{d}p_1 \cdots \mathrm{d}p_{j-1}.
$$
Here the integrals over $\R$ no longer depend on $\tau$ and we can
bound the exponential by one to obtain $t-1$ as an upper bound for
the integral with respect to $\tau$. Therefore, using~(\ref{cotaI}),
\be
\label{Nugeral}
\left|\widehat{N(u)}(w,t)\right|\leq
|\lambda|\ \frac{t-1}{1+|w|^q}\ S_0(\|u\|_L),
\ee
where $S_0$ is the
sum given by~(\ref{somaS0}). Similarly, the derivative
of~(\ref{eq:unchapeugeral}) with respect to $w$ can be bounded above
by
$$
\left(|2J(w,t)|+t-1\right)\frac{|a_j|}{(2\pi)^{j-1}}\|u\|_L^j
\int_{\R}\cdots \int_{\R} \frac{1}{1+|w-p_1|^q} \cdots
\frac{1}{1+|p_{j-1}|^q}\;\mathrm{d}p_1 \cdots \mathrm{d}p_{j-1}
$$
and using~(\ref{cotaJ}) and~(\ref{cotaI}) we conclude that
\be
\label{Nuderivgeral}
|\widehat{N(u)}'(w,t)|\leq |\lambda|
\frac{(2\sqrt{s(t)}+1)}{1+|w|^q}(t-1)S_0(\|u\|_L).
\ee
Estimates~(\ref{Nugeral}) and~(\ref{Nuderivgeral}), together with
the monotonicity of $s(t)$ and the inequality $\alpha \geq
(p+3)/(p+1)$, then imply that
\begin{eqnarray*}
\|N(u)\|_L &\leq& 2|\lambda|(\sqrt{s(L)}+1)(L-1)\|u\|_L^2 \dsum_{j
\geq \alpha}\left(\frac{C}{2\pi}\right)^{j-1}|a_{j}|\|u\|_L^{j-2}.
\end{eqnarray*}
We can now bound the sum above by its value when $\|u\|_L=\rho_0$ and
use~(\ref{ineq:u-f}) to obtain
\be
\label{eq:normaNugeral}
\|N(u)\|_L\leq C' |\lambda|\ \|f\|^2,
\ee
where
\be
\label{eq:CLnqgeral1}
C'=C'(L,q,F,c)=8\left(\sqrt{s(L)}+1\right)^3(L-1)S_1(\rho_0).
\ee
Finally, recalling that $|\lambda|\leq 1$, invoking (\ref{eq:normaf}), and
defining
$$
\varepsilon' \equiv
\dmin\left\{{C'}^{-1},[2C_q(1+\sqrt{s(L)})]^{-1}\rho_0\right\},
$$
where $C'$ is given by equation~(\ref{eq:CLnqgeral1}), we conclude
that $\|N(u)\|_L < \|f\|$ whenever $\|f\| <\varepsilon'$.
\eop

{\bf Proof of Lemma~\ref{lema:Nuvngeral}:} Consider functions $u$
and $v$ such that $\|u\|_L < \rho_0$ and $\|v\|_L< \rho_0$.
Then,
$$
[\widehat{N(u)}-\widehat{N(v)}](w,t)= \lambda \dsum_{j \geq
\alpha}a_j\int_0^{t-1}\;\mathrm{d}\tau e^{-\varphi(\tau)w^2}
[\widehat{u^j}-\widehat{v^j}](w,t-\tau) \equiv \lambda \dsum_{j\geq
\alpha}D_j,
$$
where $D_j$ can be written as
\[
D_j=\frac{a_j}{(2\pi)^{j-1}} \int_0^{t-1}\;\mathrm{d}\tau
e^{-\varphi(\tau)w^2}[(\hat u*\cdots * \hat u)- (\hat v*\cdots *
\hat v)](w,\tau).
\]
Here there are $j-1$ convolutions of $\hat u$ and $j-1$ of $\hat
v$. We add and subtract in the integrand the term $\hat v*\hat u*
\cdots * \hat u$, with $j-2$ convolutions of $\hat{u}$, to get
\begin{eqnarray*}
D_j&=& \frac{a_j}{(2\pi)^{j-1}} \int_0^{t-1}e^{-\varphi(\tau)w^2}
[(\hat u-\hat v)*\hat u* \cdots * \hat u](w,\tau)\;\mathrm{d}\tau +\\
&+& \frac{a_j}{(2\pi)^{j-1}} \int_0^{t-1}e^{-\varphi(\tau)w^2}[\hat
v*(\hat u* \cdots *\hat u- \hat v* \cdots * \hat v)](w,\tau)\;\mathrm{d}\tau.
\end{eqnarray*}
The first integral can be bounded, 
in a manner similar to what was done in the proof of
Lemma~\ref{lema:Nungeral}, by
$(t-1)(1+|w|^q)^{-1}(2\pi)^{1-j}C^{j-1}|a_j|\|u\|_L^{j-1}\|u-v\|_L$.
To estimate the second integral we rewrite it, after adding and
subtracting appropriate terms, as a sum of two integrals, one of
which can be bounded as above and the other can be split into two
other integrals. This procedure ends after $j-1$ steps, when we
obtain
$$
D_j \leq \frac{(t-1)}{1+|w|^q}\cdot
\frac{C^{j-1}|a_j|}{(2\pi)^{j-1}}\|u-v\|_L
\left(\|u\|_L^{j-1}+\|v\|_L\|u\|_L^{j-2}+ \cdots+\|v\|_L^{j-1}
\right).
$$
Note that since the norms of $u$ and $v$ in $B^{(L)}$
are less than $\rho_0$, the sum over $j \geq \alpha$ of the
right-hand side of the inequality above is convergent. In addition,
we can factor $\|u\|_L$ or $\|v\|_L$ and the remaining sum will
still be convergent. Similarly, each term of the derivative with
respect to $w$ of the difference $N(u)-N(v)$ can be written as a sum
of two integrals, which we bound using the same procedure as before.
Therefore,
\begin{eqnarray*}
\label{eq:cotadifgeral} \|N(u)-N(v)\|_L &\leq& C'' |\lambda| \|f\|
\|u-v\|_L,
\end{eqnarray*}
where \be \label{eq:CLFqgeral}
C''=C''(L,q,F,c)=4(\sqrt{s(L)}+1)^2(L-1)S_2(\rho_0). \ee Since
$|\lambda |\leq 1$, defining
\[
\varepsilon'' \equiv
\dmin\left\{(2C'')^{-1},[2C_q(1+\sqrt{s(L)})]^{-1}\rho_0\right\},
\]
where $C''$ is given by~(\ref{eq:CLFqgeral}), the lemma is proved if
we take $\|f\|<\varepsilon''$. \eop

We note for future use that, since $S_2(\rho_0)>S_1(\rho_0)$,
we may take
\be
\label{eq:C0geral}
C_0 \equiv 8(\sqrt{s(L)}+1)^{3}(L-1)S_2(\rho_0),
\ee
and it is enough to consider $\varepsilon$ in Theorem~\ref{teo:existunic} defined by
\be
\label{def:epsilon-exist}
\varepsilon \equiv
\min\left\{(2C_0)^{-1},\left[2C_q\left(\sqrt{s(L)}+1\right)\right]^{-1}\rho_0\right\}.
\ee


\section{\large{Global Existence, Uniqueness and Asymptotic Behavior}}
\label{sec:asint}

It follows from Theorem~\ref{teo:existunic} that, given $L>1$, there
exists an $\varepsilon>0$ such that the IVP~(\ref{pvi:eqc(t)F}) has
a unique solution $u$ in $B_f$ for any $f \in \B_q$ with
$\|f\|<\varepsilon$. Therefore, \be \label{eq:rg1} f_1(x)\equiv
L^{\frac{(p+1)}{2}}u\left(L^{\frac{(p+1)}{2}}x,L \right) \ee is a
well defined element of ${\cal B}_q$. The right-hand side
of~(\ref{eq:rg1}) defines an operator, $R_{L,0}$, acting on the ball
$\{f\in {\cal B}_q : \|f\|<\varepsilon\}$, which maps the initial
condition $f$ to $f_1$. We dub $R_{L,0}$ the {\em Renormalization
Group operator} associated to problem~(\ref{pvi:eqc(t)F}).

The RG operator just defined was introduced in
\cite{bib:bric-kupa-lin}; its iteration comprises the RG method for
the asymptotic analysis of solutions. The basic idea of this method
is to reduce the long-time-asymptotics problem to the analysis of a
sequence of finite-time problems obtained by iterating the RG
operator. In more detail, first consider problem~(\ref{pvi:eqc(t)F})
and, as in Section~\ref{sec:exis-uniq}, restrict the initial data so
that this problem has a unique solution. Then, apply $R_{L,0}$ to
the initial data $f$ to produce $f_1$, the initial data for a new,
{\em renormalized} IVP. It is expected that this procedure can be
iterated ad infinitum to generate a sequence of finite-time IVPs,
whose initial conditions $f_n$ are obtained by iterating the RG
operator $n$ times.

We now outline the RG method for the nonlinear problem~(\ref{pvi:eqc(t)F}).
See also \cite{bib:gas-fred-ju-leo-tp1, bib:gas-fred-ju-leo}.
Assume that the solution $u$ to IVP~(\ref{pvi:eqc(t)F}) is
globally well defined and let $L>1$ be fixed. We consider a sequence
$\{u_n\}_{n=0}^\infty$ of rescaled functions defined by \be
\label{eq:vn} u_{n}(x,t) \equiv L^{n(p+1)/2}u\left
(L^{n(p+1)/2}x,L^nt\right ), \ee with $t\in[1,L]$. A direct
calculation reveals that $u_n$ satisfies the renormalized IVP:
\be
\label{eq:unkgeral}
\cases{\partial_tu_n=c_n(t)\partial_x^2 u_n+ \lambda_n F_{n}(u_n),
& $t\in [1,L], \,\, x\in\R,$ \cr \cr u_n(x,1)=f_n(x),& $x\in\R$,}
\ee
where
$c_n(t)=L^{-np}c(L^nt)$, $\lambda_n=L^{n[p+3-\alpha(p+1)]/2}\lambda$,
$$
F_{n}(v)=\dsum_{j\geq \alpha} \left[L^{n(\alpha-j)(p+1)/2} a_j\right] v^j
$$
and
\be\label{def:fn} f_n(x) = u_n(x,1) = L^{n(p+1)/2}u\left
(L^{n(p+1)/2}x,L^n\right).
\ee
Comparing~(\ref{def:fn})
and~(\ref{eq:cotanaolin}), it becomes clear that proving the
asymptotic limit may be reduced to proving the convergence of
$\{f_n\}$, which motivates the definition of the RG operator. Let $g
\in \B_q$ and for a given $n\geq0$ assume that the
IVP~(\ref{eq:unkgeral}) with initial condition $g$ has a unique
solution $u_n$. Then, rescale $u_n(\cdot,L)$ to obtain
\begin{equation}
  \label{eq:Rndef}
  (R_{L,n}g)(x) \equiv L^{(p+1)/2}u_n\left (L^{(p+1)/2}x,L\right ),
\end{equation}
which defines the RG operator. The index $n$ in the above definition
is justified since the operator depends on the evolution equation
considered. Now if we consider IVP~(\ref{eq:unkgeral}) with initial
data $f_n$ given by~(\ref{def:fn}), then it is an immediate
consequence of these definitions that the sequence $\{f_n\}$
satisfies
\be
  \label{eq:fn_iterate}
  f_0=u(\cdot,1) {\mbox{\hspace{0.5cm} and \hspace{0.5cm}}} f_{n+1} =
  R_{L,n}f_n.
\ee
Our goal from now on is to make the above heuristic
argument rigorous. We shall prove that under hypotheses $(H1)-(H3)$,
if the initial data is sufficiently small,
problem~(\ref{eq:unkgeral}) has a unique solution for each $n$ so
that the iterative RG method can be applied to furnish the
asymptotic behavior of the solution to IVP~(\ref{pvi:eqc(t)F}).

In Lemma~\ref{lema:existuniclocal} we proceed as in the proof of
Theorem~\ref{teo:existunic} to obtain local existence and uniqueness
of solutions for each problem~(\ref{eq:unkgeral}). To state the
lemma, consider the space $B^{(L)}$ defined by~(\ref{def:spacebl})
and, if $f_n$ is the initial data of problem~(\ref{eq:unkgeral}),
define the space $B_{f_n}=\{u_n\in B^{(L)}:\|u_n-u_{f_n}\|\leq
\|f_n\|\}$ and the operator $T_n(u_n)\equiv u_{f_n} + N_n(u_n)$,
where $u_{f_n}$ is the solution of~(\ref{eq:unkgeral}) with
$\lambda_n =0$ and \be \label{eq:nukgeral} N_n(u_n)(x,t)= \lambda_n
\int_0^{t-1} \int \frac{e^{-\frac{(x-y)^2}{4[s_n(t)-s_n(t-\tau)]}}}
{\sqrt{4\pi[s_n(t)-s_n(t-\tau)]}} F_{n}(u_n(y,t-\tau))\;\mathrm{d}y\mathrm{d}\tau,
\ee where \be \label{eq:sn} s_n(t)=\int_1^t
c_n(v)\;\mathrm{d}v=\frac{t^{p+1}-1}{p+1}+r_n(t). \ee
Define also the constant $C_n$ by
\be
\label{def:Cn}
C_n \equiv
8(\sqrt{s_n(L)}+1)^3(L-1)S_2(\rho_0),
\ee
where $S_2(\rho_0)$ is given by~(\ref{somaS2}), with $z=\rho_0$.

\begin{lema}
\label{lema:existuniclocal} Given $n\in \N$ and $L>1$, there exists
an $\varepsilon_n>0$ such that if $\|f_n\|<\varepsilon_n$, then the
IVP~(\ref{eq:unkgeral}) has a unique solution $u_n(x,t)$ in
$B_{f_n}$. Furthermore, $f_{n+1}$ given by~(\ref{eq:fn_iterate}) is
a well defined element of the $\B_q$ space.
\end{lema}
\noindent {\bf Proof:}
(Notice that if $n=0$ and $f_0 \equiv f$,
Lemma~\ref{lema:existuniclocal} is just Theorem~\ref{teo:existunic}.)
We must prove that the operator $T_n$ is a
contraction in $B_{f_n}$, therefore obtaining a unique solution
$u_n$ in $B_{f_n}$. First, following closely the arguments in the
proof of Lemma~\ref{lema:Nungeral}, the constraint $L>1$ and the
definitions of $F_{n}$ and $s_n(t)$ imply that \be
\label{ineq:cotanu} \|N_n(u_n)\|_L\leq
C_nL^{n[p+3-\alpha(p+1)]/2}\|f_n\|^2 \ee and that
$$
\|N_n(u_n)-N_n(v_n)\|\leq C_n
L^{n[p+3-\alpha(p+1)]/2}\|f_n\|\|u_n-v_n\|,
$$
where $C_n$ is given by~(\ref{def:Cn}). The condition for $u_n$ to
be in the region of analyticity of $F_{n}$ is now that
$\|f_n\|<[2C_q(1+\sqrt{s_n(L)})]^{-1}\rho_0$. Since
$p+3-\alpha(p+1)<0$, defining \be \label{def:epsilonn}
\varepsilon_n \equiv
\dmin\left\{{(2C_n)}^{-1},[2C_q(1+\sqrt{s_n(L)})]^{-1}\rho_0\right\},
\ee if $\|f_n\|<\varepsilon_n$, we obtain $\|N_n(u_n)\|_L<\|f_n\|$
and $\|N_n(u_n)-N_n(v_n)\|<\frac{1}{2}\|u_n-v_n\|_L$ for all $u_n,
v_n \in \B_{f_n}$. This proves that the IVP~(\ref{eq:unkgeral}) has
a unique solution $u_n(x,t)$ in $B_{f_n}$ and, therefore,
$f_{n+1}\equiv L^{(p+1)/2}u_n\left (L^{(p+1)/2}x,L\right)$ is well
defined.\eop

We have proved that if $\|f_n\|<\varepsilon_n$, then $R_{L,n}f_n$ is
well defined. To simplify the notation, let $\nu_n(x)\equiv
N_n(u_n)(x,L)$ (cf.~(\ref{eq:nukgeral}), where $N_0(u)=N(u)$). Then,
the solution to the IVP~(\ref{eq:unkgeral}) at time $t=L$ can be
written as $u_n(x,L)=u_{f_n}(x,L)+\nu_n(x)$ and we have \be
\label{def:rg} (R_{L,n}f_n)(x)= R^0_{L}f_n(x) +
L^{(p+1)/2}\nu_n(L^{(p+1)/2}x), \ee where $R_L^0\equiv R_{L,0}^0$
and $(R^0_{L,n}f_n)(x)\equiv L^{(p+1)/2}u_{f_n}\left
(L^{(p+1)/2}x,L\right )$. We see that~(\ref{def:rg}) splits the RG
operator into two parts, which we dub the linear and the nonlinear
parts. Our analysis focus first on the linear part; the nonlinear
part is driven to zero under hypothesis $(H3)$ and, thereby,
does not contribute to the asymptotic regime, as we shall prove.

It follows from the definition of $R^0_{L,n}$ and from the integral
representation of $u_{f_n}$ that if $g \in \B_q$ is the initial data
of IVP~(\ref{eq:unkgeral}) with $\lambda_n=0$, then the Fourier
Transform of $R^0_{L,n}g$ is given by \be
  \label{eq:rgn}
 \F({R^0_{L,n}g})(w)= \widehat{g}(L^{-(p+1)/2}w
)e^{-w^2s_n(L)/L^{p+1}}, \ee where $s_n$ is defined
in~(\ref{eq:sn}). Applying equation~(\ref{eq:rgn}) inductively and
using that $s_0(L^n)+L^{n(p+1)}s_n(L)=s_0(L^{n+1})$ for all
$n=1,2,\dots$, it is easy to prove that the linear RG operator has
the semi-group property. Also, if $g=f_p^*$, it follows from
equation~(\ref{eq:rgn}) and definition~(\ref{eq:sn}) with $n=0$ that
\be \label{eq:rgfp} \F({R^0_{L}f_p^*})(w) =  e^{-w^2/(p+1)}
e^{-w^2r(L)/L^{p+1}}, \ee where $r(L)\equiv r_0(L)$  (see
(\ref{eq:sn}) for the definition of $r_0(L)$) and we have used that
$\hat{f_p^*}(w)=e^{-w^2/(p+1)}$. Taking the inverse Fourier
Transform on both sides of equation~(\ref{eq:rgfp}) we conclude that
$f_p^*(x)$ is a fixed point of the RG operator if and only if $r(L)
\equiv 0$, which is valid if we take $c(t) = t^p$. The next lemma
shows that if $r(t)\not \equiv 0$, $f_p^*$ is still the long-time
asymptotic limit of $R_{L}^0f_p^*$ and, therefore, it will be an
{\it asymptotic fixed point} of the linear RG operator.
\begin{lema}
\label{teo:asym-fixe-poin} There is a constant $M=M(p,q)$ and an
$n_0 \in \N$ such that \be \label{ineq:converg}
\|R_{L^n}^{0}f_p^*-f_p^*\| \leq M|L^{-n(p+1)}r(L^{n})| \ee for all
$n>n_0$. In particular, $\|R^0_{L^{n}}f_p^* - f_p^*\|
{\longrightarrow}0 \,\,\ \mbox{as}\,\,\, n\to\infty.$
\end{lema}
\noindent {\bf Proof:} Since $r(t)=o(t^{p+1})$, from the semi-group
property of $R_{L}^0$ and from equation~(\ref{eq:rgfp}), with $L$
replaced by $L^{n}$, we have pointwise convergence in Fourier space.
To prove the theorem, we have to show convergence in $\B_q$. This
will be done at the same time that we estimate the rate of
convergence. From~(\ref{eq:rgfp}),
$$
 \left|[\F(R_{L^n}^{0}f_p^*-f_p^*)](w)\right|
 \leq  w^2\left|\frac{r(L^{n})}{L^{n(p+1)}}\right|
e^{-w^2\left[\frac{1}{p+1}-\left|\frac{r\left(L^{n}\right)}{L^{n(p+1)}}\right|
\right]}
$$
and
$$
\left|[\F{(R^0_{L^n}f_p^*-f_p^*)}]'(w)\right|\leq
2\left[\frac{|w|^3}{p+1}+|w|\right]
\left|\frac{r(L^{n})}{L^{n(p+1)}}\right|
e^{-w^2\left[\frac{1}{p+1}-\left|\frac{r\left(L^{n}\right)}{L^{n(p+1)}}\right|
\right]}.
$$
Since $r(t) = o(t^{p+1})$, there exists an $n_0>0$ such that
$\left|r(L^{n})L^{-n(p+1)}\right|<[2(p+1)]^{-1}$ for all $n>n_0$.
Multiplying the inequalities above by $(1 + |w|^q)$ and defining
$M\equiv
\max_w(1+|w|^q)e^{-[2(p+1)]^{-1}w^2}[2|w|+w^2+2|w|^3/(p+1)]$ we
obtain (\ref{ineq:converg}). Letting $n \to \infty$ finishes the
proof. \eop

In the next lemma we prove that if $L$ is sufficiently large, then
the linear RG operator $R_{L,n}^0$ is a contraction in the space of
functions $g \in \B_q$ such that $\hat g(0)=0$. This result will be
used to obtain the estimates of the Renormalization
Lemma~\ref{lema1geral}, which will allow us to prove
Theorem~\ref{teo:compassgeral}.
\begin{lema}[Contraction Lemma]
\label{lema:contracao}
\noindent There exist constants $L_1>1$ and $C>0$ such
that the inequality
\be
\label{ineq:contracao}
 \left\|R_{L,n}^0g\right\| \leq L^{-(p+1)/2}C \|g\|
\quad (n=0,1,\dots)
\ee
holds if $L>L_1$ and $g \in \B_{q}$ satisfies $\hat{g}(0) = 0$.
\end{lema}
{\bf Proof:} We bound~(\ref{eq:rgn}) and its derivative by
$$
e^{-w^2\frac{s_n(L)}{L^{p+1}}}\left[
 \left|\hat{g}\left(\frac{w}{L^{\frac{p+1}{2}}}\right)\right|+
2|w|\left|\frac{s_n(L)}{L^{p+1}}\right| \left|\hat{g}
  \left(\frac{w}{L^{\frac{p+1}{2}}}\right)\right|+
L^{-(p+1)/2}\left|\hat{g}'\left(\frac{w}{L^{\frac{p+1}{2}}}\right)\right|
 \right].
$$
From~(\ref{eq:sn}), it is easy to see that
$r_n(t)=L^{-n(p+1)}[r(L^nt)-r(L^n)]$. Now, since $r(t) =
o(t^{p+1})$, we conclude that there exists an $L_0>1$ such that, if
$L>L_0$, then $|L^{-(p+1)}r_n(L)|<[2(p+1)]^{-1}$ for all $n$. This
together with the definition  of $s_n$ furnishes \be
\label{eq:cotasn} \frac{1}{6(p+1)} \leq \frac{s_n(L)}{L^{p+1}}  \leq
\frac{3}{2(p+1)} \quad (n=0,1,2,\dots) \ee if $L>L_1\equiv
\max\{L_0,\sqrt[p+1]{3}\}$. Since $g \in \B_q$, from the definition
of the $\B_q$ norm,
$\left|\hat{g}'\left(L^{-(p+1)/2}w\right)\right|\leq \|g\|$ and
since $\hat{g}(0) = 0$ it follows that
$\left|\hat{g}\left(L^{-(p+1)/2}w\right)\right|\leq
L^{-(p+1)/2}|w|\|g\|$. Using these bounds, if $L>L_1$ we obtain, for
all $n=0,1,2,\dots,$
$$
|(\widehat{R^0_{L,n}g})(w)|+ |(\widehat{R^0_{L,n}g})'(w)|
  \leq
  \left(1+|w| + 3(p+1)^{-1}|w|^2\right)
e^{\frac{-w^2}{6(p+1)}} L^{-(p+1)/2}\|g\|.
$$
Defining $C\equiv \sup_w\left(1+|w| +
3(p+1)^{-1}|w|^2\right)(1+|w|^q) e^{\frac{-w^2}{6(p+1)}}$ finishes
the proof. \eop

The RG analysis involves decomposing the initial condition into two
factors: one in the direction of the asymptotic fixed point of the
RG operator and the other in a direction which is {\it irrelevant in
the RG sense}. That is, when the RG operator is applied to the
initial data, the irrelevant component is contracted for large $L$.

In Lemma~\ref{lema1geral} we decompose $f_n$ given
by~(\ref{eq:fn_iterate}) and obtain estimates needed later to prove
Theorem~\ref{teo:compassgeral}. We will first assume that, given
$k\in \N$,  $f_n$ is well defined for all $n=0,1,\dots,k$, which is
guaranteed if $\|f_n\|<\varepsilon_n$ for all $n=0,1,\dots,k-1$ (cf.
Lemma~\ref{lema:existuniclocal}). Later, in
Lemma~\ref{teo:estimativageral} we prove that if $f_0$ is small
enough, then the sequence $\{f_n\}$ given by~(\ref{eq:fn_iterate})
is well defined. Furthermore, notice that from the definition of the
$\B_q$ norm and some estimates used in the proof of
Lemma~\ref{lema:contracao}, we obtain constants $C_{p,q}$ and
$K_{p,q}$ such that $\|f_p^*\|\leq C_{p,q}$ and, if $L>L_1$, \be
\label{eq:cotafp} \|R^0_{L^n}f_p^*\|\leq K_{p,q}, \,\,\ \forall
n=1,2,\dots \ee For the next lemmas, we will always refer to the
constants $C_{p,q}$ and $K_{p,q}$ and to $C$ and $L_1$ given in
Lemma~\ref{lema:contracao}.

\begin{lema}[Renormalization Lemma]
\label{lema1geral}
Given $L>L_1$, suppose $f_n$ is defined for $n=0,1,\dots,k+1$ as specified
in $(\ref{eq:fn_iterate})$.  We can then write
\be
\label{decompfk} f_0=A_0f_p^*+g_0, \quad
f_{n+1}=A_{n+1}R^0_{L^{n+1}}f_p^*+g_{n+1} \quad (n=0,1, \dots, k)
\ee
in terms of functions $g_n \in \B_q$, $\hat g_{n}(0)=0$
and constants $A_n$, $K$ which satisfy
\be
\label{ineq:gn} \|g_{n+1}\| \leq C L^{-(p+1)/2}\|g_n\|+ K
L^{n[p+3-\alpha(p+1)]/2}\|f_n\|^2
\ee
and
\be
\label{item1geral}
|A_{n+1}-A_n|\leq C_n L^{n[p+3-\alpha(p+1)]/2}\|f_n\|^2,
\ee
with $C_n$ given by~(\ref{def:Cn}).
\end{lema}

{\bf Proof:} We first prove (\ref{decompfk}) inductively. Define
$g_0$ by $f_0=A_0f_p^*+g_0$, with $A_0=\hat f_0(0)$ and since
$\widehat{f_p^*}(0)=1$, we have $\hat g_0(0)=0$. By hypothesis,
$f_1$ is well defined by $R_{L,0}f_0$ and using
representation~(\ref{def:rg}) and the decomposition above for $f_0$
we can write $f_1=A_1R^0_Lf_p^*+g_1$, where $A_1=A_0+\hat{\nu_0}(0)$
and
$g_1(x)=R^0_Lg_0(x)+L^{(p+1)/2}\nu_0(L^{(p+1)/2}x)-\hat{\nu_0}(0)R^0_Lf_p^*(x)$.
It follows from the definition of $R_L^0$ that $\F(R^0_Lg_0)(0)=0$
and $\F(R^0_Lf_p^*)(0)=1$ and, therefore, $\hat g_1(0)=0$, which
proves~(\ref{decompfk}) for $n=0$. Now suppose~(\ref{decompfk})
holds for $n=0,\dots,j-1$, where $j\leq k$. We will prove that it
holds also for $n=j$. Using~(\ref{decompfk}) with $n=j-1$,
representation~(\ref{def:rg}) and the semi-group property of the
linear RG operator we obtain \be \label{eq:fk+1geral}
f_{j+1}(x)=A_jR^0_{L^{j+1}}f_p^*(x)+R^0_{L}g_j(x)+L^{(p+1)/2}\nu_j(L^{(p+1)/2}x).
\ee Defining \be \label{eq:Ak+1geral} A_{j+1}\equiv
A_j+\hat{\nu_j}(0) \ee and \be \label{eq:gk+1geral} g_{j+1}(x)
\equiv R^0_{L}g_j(x)+L^{(p+1)/2}\nu_j(L^{(p+1)/2}x)-
\hat{\nu_j}(0)R^0_{L^{j+1}}f_p^*(x), \ee we can
write~(\ref{eq:fk+1geral}) as
$f_{j+1}=A_{j+1}R^0_{L^{j+1}}f_p^*+g_{j+1}.$ From the induction
hypothesis, $\hat g_j(0)=0$ and therefore, from
definition~(\ref{eq:gk+1geral}), since the Fourier Transforms of
${R^0_{L}g_j}$ and $R^0_{L^{j+1}}f_p^*$ at the origin are equal,
respectively to $\hat g_j(0)$ and $\hat{f_p^*}(0)$, we obtain $\hat
g_{j+1}(0)=0$, which proves~(\ref{decompfk}) for $n=0,1,\dots,k$.

Recalling that $\nu_n(x) \equiv N_n(u)(x,L)$ and since
estimate~(\ref{ineq:cotanu}) holds for all $n$, using
definition~(\ref{eq:Ak+1geral}) we obtain~(\ref{item1geral}) for
$n=0,1,\dots,k$. After a calculation similar to the one in the proof
of Lemma~\ref{lema:Nungeral}, we obtain
$\|L^{(p+1)/2}\nu_n(L^{(p+1)/2}\cdot)\| \leq
L^{(p+1)q/2}C_nL^{n[p+3-\alpha(p+1)]/2}\|f_n\|^2$ and
using~(\ref{eq:cotafp}),
$\|L^{(p+1)/2}\nu_n(L^{(p+1)/2}\cdot)-\hat{\nu_n}(0)R^0_{L^{n+1}}f_p^*\|\leq
(L^{(p+1)q/2}+K_{p,q})C_n L^{n[p+3-\alpha(p+1)]/2}\|f_n\|^2$. Since
$\hat g_n(0)=0$ and $L>L_1$, from definition~(\ref{eq:gk+1geral})
and Lemma~\ref{lema:contracao}, we obtain
\be
\label{item2geral}
\|g_{n+1}\| \leq C L^{-(p+1)/2}\|g_n\|+
\left(L^{(p+1)q/2}+K_{p,q}\right)C_n
L^{n[p+3-\alpha(p+1)]/2}\|f_n\|^2
\ee
for all $n=0,1,\dots,k$. Now it follows from~(\ref{eq:cotasn}) that
the constants $C_n$ are uniformly bounded. In fact, defining
\be
\label{def:K}
K\equiv 8(L-1)\left(\sqrt{\frac{3L^{p+1}}{2(p+1)}}+1 \right)^3
\left(L^{(p+1)q/2}+K_{p,q}\right)S_2(\rho_0),
\ee
then $C_n\leq K$ for all $n$ and therefore we obtain inequality~(\ref{ineq:gn}),
which ends the proof.

\eop

The estimates obtained in
Lemma~\ref{lema1geral} are used to prove
Theorem~\ref{teo:compassgeral} in the following
way:~(\ref{item1geral}) guarantees that the sequence $(A_n)$ is
convergent and~(\ref{ineq:gn}) is used to prove that the
component $g_n$ gets smaller as we increase $n$. This is so because
of our definition of $\alpha$ or, in other words, because the
nonlinear perturbation $F$ of problem~(\ref{pvi:eqc(t)F}) is
irrelevant. Before we apply Lemma~\ref{lema1geral}, we have to prove
that the initial data of each problem~(\ref{eq:unkgeral}) is small
enough and to do that we will define a recursive sequence $(G_n)$
such that, for all $n$, $\|f_n\|\leq G_n\|f_0\|$. In
Lemma~\ref{lema:Gn} we prove that, under a certain condition, this
sequence is bounded. Given $\delta \in (0,1)$, let \be
\label{eq:ldeltageral} L_{\delta}\equiv
\max\{L_1,[2C(1+C_{p,q})]^{2/(p+\delta)}\} \ee and for
$L>L_{\delta}$, define \be \label{def:G} G \equiv
1+K_{p,q}\sum_{j=0}^{\oo}L^{j(\delta-1)/2}<\oo, \ee $G_1\equiv
L^{(\delta-1)/2}+ K_{p,q}(1+C_0\|f\|)$ and $G_{n+1}$, for
$n=1,2,3,\dots,$ by the relation:
$$
G_{n+1}\equiv L^{(\delta-1)(n+1)/2}+ K_{p,q}\left(1+C_0\|f\|+
\dsum_{j=1}^{n}C_j G_j^{2}L^{j[p+3-\alpha(p+1)]/2}\|f\|\right),
$$
where each $C_j$, with $j=0,1,2,\dots$, is given by
equation~(\ref{def:Cn}), with $n=j$.

\begin{lema}
\label{lema:Gn} Let $\delta \in (0,1)$ be such that
$\delta-1>p+3-\alpha(p+1)$ and let $L >L_{\delta}$, where
$L_{\delta}$ is given by~(\ref{eq:ldeltageral}). Also, let $K$ and
$G$ be given by~(\ref{def:K}) and~(\ref{def:G}), respectively, and
suppose that $f$ satisfies \be \label{eq:hipED2f0geral}
KG^2\|f\|<\frac{1}{2L^{(1-\delta)/2}}. \ee Then $G_{n+1}< G$ for all
$n=0,1,2,\dots$ .
\end{lema}

\noindent {\bf Proof:} Since $L>1$ and $G>1$, it is straightforward
from the fact that $C_0 \leq K$ and from
condition~(\ref{eq:hipED2f0geral}) that $G_1<G$. Now, suppose
$G_{n+1}< G$, $\forall n=1,2,\dots,k-1$. From the definition of
$G_{n+1}$, using the induction hypothesis and since $C_n \leq K$,
$\forall n$,
$$
G_{k+1}\leq L^{(\delta-1)(k+1)/2}+K_{p,q}\left(
1+K\|f\|+KG^2\|f\|\dsum_{j=1}^{k}L^{j[p+3-\alpha(p+1)]/2}\right).
$$
Now, from~(\ref{eq:hipED2f0geral}) and since $L>1$ and
$\delta-1>p+3-\alpha(p+1)$, we obtain $G_{k+1}\leq 1+
K_{p,q}\left(1+L^{(\delta-1)/2}+\cdots+
L^{(k+1)(\delta-1)/2}\right)<G$, which completes the proof. \eop

In Lemma~\ref{teo:estimativageral} we will obtain estimates for the
rescaled solutions to IVP~(\ref{eq:unkgeral}). In fact, we will
define $\overline\varepsilon>0$ such that, if the initial data $f$
of problem~(\ref{pvi:eqc(t)F}) is in the ball of radius
$\overline\varepsilon$, then there is a unique global solution to
IVP~(\ref{pvi:eqc(t)F}). Furthermore, we will prove that, under
certain hypotheses, the component $g_n$ of the initial data $f_n$
goes to zero when $n \to \infty$. This fact will be used to obtain
the asymptotic behavior in Theorem~\ref{teo:asint1}. Before stating
the lemma, we notice that from~(\ref{eq:cotasn}), if $\varepsilon_n$
is given by~(\ref{def:epsilonn}) and
\be
\label{def:sigma}
\sigma
\equiv \dmin\left\{{(2K)}^{-1},
\left[2C_q\left(1+\sqrt{\frac{3L^{p+1}}{2(p+1)}}\right)\right]^{-1}\rho_0\right\},
\ee
then $\sigma<\varepsilon_n$ for all $n$. In the next Lemma we
will refer to $K$, $L_{\delta}$, $G$ and $\sigma$ given,
respectively,
by~(\ref{def:K}),~(\ref{eq:ldeltageral}),~(\ref{def:G})
and~(\ref{def:sigma}).
\begin{lema}
\label{teo:estimativageral} Let $L>L_{\delta}$ and $\delta \in
(0,1)$ such that $\delta-1>p+3-\alpha(p+1)$. Then, there is
$\overline{\varepsilon}>0$ such that, if
$\|f_0\|<\overline\varepsilon$, $f_{n}$ given
by~(\ref{eq:fn_iterate}) is well defined for all $n\geq 1$, \be
\label{eq:cotafkgeral} \|f_n\|\leq G_n\|f_0\| \ee and if $g_n$ is
given by the decomposition~(\ref{decompfk}), then, \be
\label{eq:cotagkgeral} \|g_n\|\leq
\frac{\|f_0\|}{L^{n(1-\delta)/2}}. \ee
\end{lema}

\noindent {\bf Proof:} The proof is by induction in $n$. First,
define \be \label{eq:epsilonbarrageral} \overline{\varepsilon}
\equiv \dmin \left\{\sigma G^{-1},
[2KG^2L^{(1-\delta)/2}]^{-1}\right\}. \ee Since $G>1$, we have
$\|f_0\|<\sigma<\varepsilon_0$ and, from
Lemmas~\ref{lema:existuniclocal} and~\ref{lema1geral}, $f_1$ is well
defined by $f_1=A_1R^0_{L}f_p^*+g_1$ and $g_1$
satisfies~(\ref{ineq:gn}) with $k=0$. Therefore, since
$f_0=A_0f_p^*+g_0$, we obtain $\|g_1\|\leq
[C(1+C_{p,q})L^{-(p+1)/2}+K\|f_0\|]\|f_0\|$. Since $L>L_{\delta}$,
then, $2C(1+C_{p,q})L^{-(p+1)/2}<L^{(\delta-1)/2}$ and since $G>1$
and $\|f_0\|<\overline\varepsilon$, then
$2K\|f_0\|<L^{(\delta-1)/2}$. Therefore, $\|g_1\| \leq
L^{(\delta-1)/2}\|f_0\|$. Now, using decomposition~(\ref{decompfk})
with $k=0$ and the bound~(\ref{eq:cotafp}),
$$
\|f_1\| \leq \left[(1+C_0\|f_0\|)K_{p,q}+
L^{(\delta-1)/2}\right]\|f_0\| = G_1\|f_0\|,
$$
which proves the Theorem for $n=1$. Now suppose there exists $k>1$
such that~(\ref{eq:cotafkgeral}) and~(\ref{eq:cotagkgeral}) hold for
all $n=1,2,\dots,k$. We will prove that these estimates hold also
for $n=k+1$. From the induction hypothesis and Lemma~\ref{lema:Gn},
since $\|f_0\|<\overline\varepsilon$, we have $\|f_n\|\leq
G\|f_0\|<\varepsilon_n$, $\forall n=1,2,\dots,k.$ Therefore, we can
apply Lemma~\ref{lema1geral} to obtain estimate~(\ref{ineq:gn}) with
$n=k$. Then, using~(\ref{eq:cotafkgeral}) and~(\ref{eq:cotagkgeral})
with $n=k$, we get:
$$
\|g_{k+1}\|\leq
L^{(\delta-1)(k+1)/2}\left[\frac{C}{L^{(p+\delta)/2}}+
\frac{L^{k[p+3-\alpha(p+1)]/2}}{L^{(\delta-1)(k+1)/2}} K
G_k^2\|f_0\|\right]\|f_0\|.
$$
Since $C_{p,q}>0$ and $L>L_{\delta}$, then $C L^{-(p+\delta)/2}<1/2$
and since $\|f_0\|<\overline\varepsilon$, using Lemma~\ref{lema:Gn}
we obtain~(\ref{eq:cotagkgeral}) with $n=k+1$. By
Lemma~\ref{lema1geral}, $f_{k+1}$ is well defined and can be
represented by~(\ref{decompfk}). Therefore, using the triangle
inequality and~(\ref{eq:cotafp}),
\begin{eqnarray*}
\|f_{k+1}\|&\leq& \frac{\|f_0\|}{L^{(1-\delta)(k+1)/2}}+
K_{p,q}\left(|A_0|+\sum_{j=0}^{k}|A_{j+1}-A_j| \right).
\end{eqnarray*}
Now, since $|A_0|\leq \|f_0\|$ and $C_n \leq K$, for all $n$,
applying estimates~(\ref{item1geral}) and~(\ref{eq:cotafkgeral}) for
$n=0,1,2,\dots,k$ and using Lemma~\ref{lema:Gn}, we
obtain~(\ref{eq:cotafkgeral}) with $n=k+1$. In particular,
$\|f_{k+1}\|<G\|f_0\|<\varepsilon_{k+1}$, which ends the proof. \eop

We have proved that for $\|f_0\|<\overline\varepsilon$ each
IVP~(\ref{eq:unkgeral}) has a unique solution $u_n$ in $B_{f_n}$.
To finish the proof we only need to concatenate these solutions
to obtain a unique global solution to IVP~(\ref{pvi:eqc(t)F}).

We first extend the
definition~(\ref{def:spacebl}) of the $B^{(L)}$ space by considering
the space
$$
B^{(L^{n+1})}\equiv \left\{u:\R \times [L^n,L^{n+1}] \rightarrow \R \,\ |\,\
\|u\|_{L^{n+1}} < \oo
\right\}
$$
with the norm $\|u\|_{L^{n+1}}=\sup_{t \in
[L^{n},L^{n+1}]}\|u(\cdot,t)\|$.

Now define $\{h_n\}$ by
$$
h_0\equiv f {\mbox{\hspace{0.5cm} and \hspace{0.5cm}}} h_{n+1}\equiv
L^{-n(p+1)/2}u_n \left(L^{-n(p+1)/2}x,L\right)
$$
and let $u_{h_n}$ be the solution to IVP~(\ref{eq:unkgeral}) with
$\lambda_n=0$ and initial condition $h_n$. Finally, define
\be
\label{def:spaceb}
B\equiv \left\{u\in
B^{(\infty)}:\|u-u_{h_n}\|_{L^{n+1}}\leq \|h_n\|, \forall n\right\},
\ee
where (abusing notation) $\|\cdot\|_{L^{n+1}}$ denotes the
seminorm induced by the obvious projection from $B^{(\oo)}$ onto
$B^{(L^{n+1})}$.

\begin{cor}
\label{cor:existglob} Under the hypotheses of
Lemma~\ref{teo:estimativageral} the IVP~(\ref{pvi:eqc(t)F}) has a
unique solution $u \in B$. In this case, the RG transformation has
the \emph{``semi-group property''}: $R_{L^n,0}f = R_{L,n-1}\circ \cdots \circ
R_{L,1} \circ R_{L,0}f$ for all $n\geq 1$.
\end{cor}

\noindent {\bf Proof:} From Lemma~\ref{teo:estimativageral}, since
$\|f_0\|<\overline\varepsilon$, then $\|f_n\|<\varepsilon_n$ for all
$n$ and using Lemma~\ref{lema:existuniclocal}, we obtain the
existence and uniqueness of solutions to
problems~(\ref{eq:unkgeral}) in $B_{f_n}$. Now define $ u(x,t)
\equiv L^{-n(p+1)/2}u_n \left(L^{-n(p+1)/2}x,L^{-n}t\right), \,\,\
\forall t \in [L^{n},L^{n+1}] $ and take $y=L^{-n(p+1)/2}x$ and
$\tau=L^{-n}t$. Since $u_n(y,\tau)$ is the unique solution to
IVP~(\ref{eq:unkgeral}) in $B_{f_n}$, then $u(x,t)$ is the unique
solution to IVP~(\ref{pvi:eqc(t)F}) in $B$. To prove the semi-group
property, it is enough to apply Lemma~\ref{lema:existuniclocal}
and~(\ref{eq:Rndef}), inductively. \eop

The previous results are concatenated in Theorem~\ref{teo:asint1}.
\begin{teo}
\label{teo:asint1} Under the hypotheses of
Lemma~\ref{teo:estimativageral}, there is a constant $A$ such that
$$
\dlim_{n\to\oo}\left\|L^{n(p+1)/2}u(L^{n(p+1)/2}
\cdot,L^{n})-Af_p^*(\cdot)\right\|=0.
$$
\end{teo}

\noindent {\bf Proof:} Since $\|f_0\|<\overline\varepsilon$, it
follows from Corollary~\ref{cor:existglob} that the
IVP~(\ref{pvi:eqc(t)F}) has a unique solution $u\in B$. It follows
from the semi-group property and Lemma~\ref{lema1geral} that
$f_n=A_nR^0_{L^{n}}f_p^*+g_n=L^{n(p+1)/2}u(L^{n(p+1)/2}x,L^{n})$.
Therefore, estimate~(\ref{eq:cotagkgeral}) can be written as
$\|f_n-A_nR^0_{L^{n}}f_p^*\|\leq L^{n(\delta-1)/2}\|f_0\|$. Since
$C_n \leq K$ and $\|f_0\|<\overline\varepsilon$, using
Lemma~\ref{lema:Gn} and estimates~(\ref{item1geral})
and~(\ref{eq:cotafkgeral}), we obtain, for all $n=1,2,3,\dots$,
$$
|A_{n+1}-A_n|<
\frac{L^{n[p+3-\alpha(p+1)]/2}}{2L^{(1-\delta)/2}}\|f_0\|
$$
and therefore, $(A_n)$ is a Cauchy sequence in $\R$. Let $A \in \R$
be the limit of this sequence. Using the triangle inequality,
$\|L^{n(p+1)/2}u(L^{n(p+1)/2}x,L^{n})-Af_p^*\| \leq
\|L^{n(p+1)/2}u(L^{n(p+1)/2}x,L^{n})-A_nR_{L^n}^0f_p^*\|+|A|\|R_{L^n}^0f_p^*-f_p^*\|+
 |A_n-A|\|R_{L^n}^0f_p^*\|$, which, from
Lemma~\ref{teo:asym-fixe-poin} and~(\ref{eq:cotafp}) can be upper
bounded by \be \label{finalbound}
\frac{\|f_0\|}{L^{n(1-\delta)/2}}+|A|M\left|\frac{r(L^{n})}{L^{n(p+1)}}\right|
+ \frac{L^{n[p+3-\alpha(p+1)]/2}}{2L_{\delta}^{(1-\delta)/2}
\left(1-L_{\delta}^{[p+3-\alpha(p+1)]/2}\right)}K_{p,q}\|f_0\|. \ee
Then it is enough to take the limit when $n \to \oo$. \eop

Theorem~\ref{teo:compassgeral} now follows from estimate
(\ref{finalbound}) as we explain below.

\noindent {\bf Proof of Theorem~\ref{teo:compassgeral}:} We have
proved that~(\ref{eq:cotanaolin}) is valid when the initial data
$f\equiv f_0$ is sufficiently small and $t=L^{n}$ ($n=1,2,\dots$),
for $L>L_{\delta}$. In fact, it follows from (\ref{finalbound})
that, if $t=L^{n}$, then
$\|\sqrt{t^{p+1}}u(\sqrt{t^{p+1}}x,t)-Af_p^*\|\leq t^{(\delta-1)/2}
\|f_0\| + |A|M\left|t^{-(p+1)}r(t)\right| + t^{[p+3-\alpha(p+1)]/2}
L_{\delta}^{(\delta-1)/2}
\left[2\left(1-L_{\delta}^{[p+3-\alpha(p+1)]/2}\right)\right]^{-1}K_{p,q}\|f_0\|$.
We can extend this bound to $t = \tau L^n$, with $\tau \in [1,L]$
and $L>L_{\delta}$ by replacing everywhere $L$ by $\tau^{1/n} L$.
Therefore, since the constants in (\ref{finalbound}) do not depend
on the particular value of $L>L_\delta$ considered, the inequality
above holds for all $t>L_{\delta}$. Taking the limit $t \to \infty$
finishes the proof. \eop

To conclude, we remark that we do not have an explicit expression for the
limit $A$ of the sequence $A_n$. However, it should be clear that
\[
A = \lim_{t\to\infty}\int u(x,t)\;\mathrm{d}x.
\]
Also, we point out that the RG approach has been used to obtain
higher order corrections to the asymptotic behavior, see
\cite{bib:bona1,bib:bric-kupa}).



{\em E-mail address:} gbraga@mat.ufmg.br\\
{\em E-mail address:} furtado@uwyo.edu\\
{\em E-mail address:} jmoreira@mat.ufmg.br\\
{\em E-mail address:} leorolla@impa.br



\parskip 0pt
\baselineskip = 18pt

\begin{thebibliography}{10}

\bibitem{bib:aronson-vazquez}
D.~G. Aronson and J.~L. Vazquez.
\newblock Calculation of anomalous exponents in nonlinear diffusion.
\newblock {\em Physical Review Letters}, 17:348--351, 1994.

\bibitem{bib:avell}
M.~Avellaneda.
\newblock Homogenization and renormalization: the mathematics of multi-scale
  random media and turbulent diffusion.
\newblock In {\em Dynamical Systems and Probabilistic Methods in Partial
  Differential Equations}, volume~31 of {\em Lectures in Applied Mathematics},
  pages 251--268, 1994.

\bibitem{bib:barenblatt-book}
G.~I. Barenblatt.
\newblock {\em Scaling, self-similarity and intermediate asymptotics}.
\newblock Cambridge University Press, Cambridge, 2 edition, 1996.

\bibitem{bib:bogoliubov}
N.~N. Bogoliubov and D.~V. Shirkov.
\newblock {\em Introduction to the Theory of Quantized Fields}.
\newblock Interscience Publishers Ltd, New York, 1959.

\bibitem{bib:bona1}
J.~L. Bona, K.~S. Promislow, and C.~E. Wayne.
\newblock Higher order asymptotics of decaying solutions of some nonlinear,
  dispersive, dissipative wave equations.
\newblock {\em Nonlinearity}, 8:1179--1206, 1995.

\bibitem{bib:gas-fred-ju-leo-tp1}
G.~A. Braga, F.~Furtado, J.~M. Moreira, and L.~T. Rolla.
\newblock Renormalization group analysis of nonlinear diffusion equations with
  time dependent diffusion coefficients: Numerical results.
\newblock Paper in preparation.

\bibitem{bib:gas-fred-ju-leo}
G.~A. Braga, F.~Furtado, J.~M. Moreira, and L.~T. Rolla.
\newblock Renormalization group analysis of nonlinear diffusion equations with
  periodic coefficients.
\newblock {\em Multiscale Modeling and Simulation}, 1(4):630--644, 2003.

\bibitem{bib:bric-kupa}
J.~Bricmont and A.~Kupiainen.
\newblock Renormalizing partial differential equations.
\newblock {\em Constructive Physics}, 446:83--115, 1995.

\bibitem{bib:bric-kupa-lin}
J.~Bricmont, A.~Kupiainen, and G.~Lin.
\newblock Renormalization group and asymptotics of solutions of nonlinear
  parabolic equations.
\newblock {\em Communications in Pure and Applied Mathematics}, 47:893--922,
  1994.

\bibitem{bib:bric-kup-xin}
J.~Bricmont, A.~Kupiainen, and J.~Xin.
\newblock Global large time self-similarity of a thermal-diffusive combustion
  system with critical nonlinearity.
\newblock {\em J. Diff. Eq.}, 130:19--35, 1996.

\bibitem{bib:cag}
G.~Caginalp.
\newblock Renormalization and scaling methods for nonlinear parabolic systems.
\newblock {\em Nonlinearity}, 10:1217--1229, 1997.

\bibitem{bib:gold3}
L.-Y. Chen, N.~Goldenfeld, and Y.~Oono.
\newblock Renormalization group and singular perturbations: multiple scales,
  boundary layers, and reductive perturbation theory.
\newblock {\em Physical Review E}, 54:376--394, 1996.

\bibitem{bib:dagan87}
G.~Dagan.
\newblock Theory of solute transport by groundwater.
\newblock {\em Ann. Rev. Fluid Mech.}, 19:183--215, 1987.

\bibitem{bib:davis91}
R.~E. Davis.
\newblock Lagrangian ocean studies.
\newblock {\em Ann. Rev. Fluid Mech.}, 23:43, 1991.

\bibitem{bib:gell-low}
M.~Gell-Mann and F.~E. Low.
\newblock Quantum electrodynamics at small distances.
\newblock {\em Phys. Rev.}, 95:1300--1312, 1954.

\bibitem{bib:glimm}
J.~Glimm, W.~B. Lindquist, F.~Pereira, and Q.~Zhang.
\newblock A theory of macrodispersion for the scale-up problem.
\newblock {\em Transport in Porous Media}, 13:97--122, 1993.

\bibitem{bib:gold-book}
N.~Goldenfeld.
\newblock {\em Lectures on Phase Transition and the Renormalization Group}.
\newblock Addison-Wesley, Reading, 1992.

\bibitem{bib:gold2}
N.~Goldenfeld, O.~Martin, and Y.~Oono.
\newblock Asymptotics of partial differential equations and the renormalization
  group.
\newblock In S.~Tanveer, editor, {\em Proc. NATO Advanced Research Workshop on
  Asymptotics Beyond All Orders}, New York, 1992. Plenum.

\bibitem{bib:gold1}
N.~Goldenfeld, O.~Martin, Y.~Oono, and F.~Liu.
\newblock Anomalous dimensions and the renormalization group in a nonlinear
  diffusion process.
\newblock {\em Phys. Rev. Lett.}, 64:1361--1364, 1990.

\bibitem{bib:lake86}
L.~Lake and H.~Carroll, editors.
\newblock {\em Reservoir Characterization}.
\newblock Academic Press, New York, 1986.

\bibitem{bib:li-qi}
Y.~Li and Y.~W. Qi.
\newblock The global dynamics of isothermal chemical systems with critical
  nonlinearity.
\newblock {\em Nonlinearity}, 16:1057--1074, 2003.

\bibitem{bib:majda}
A.~J. Majda and P.~R. Kramer.
\newblock Simplified models for turbulent diffusion: theory, numerical
  modelling, and physical phenomena.
\newblock {\em Physics Reports}, 314:237--574, 1999.

\bibitem{bib:merd}
H.~Merdan and G.~Caginalp.
\newblock Decay of solutions to nonlinear parabolic equations: Renormalization
  and rigorous results.
\newblock {\em Discrete and Continuous Dynamical Systems B}, 3:565--588, 2003.

\bibitem{bib:moffatt83}
H.~K. Moffatt.
\newblock Transport effects associated with turbulence.
\newblock {\em Rep. Prog. Phys}, 46:621, 1983.

\bibitem{bib:oono}
Y.~Oono.
\newblock Renormalization and asymptotics.
\newblock {\em Int. Journal of Modern Physics B}, 14:1327--1361, 2000.

\bibitem{bib:vazquez-ima}
J.~L. Vazquez.
\newblock Asymptotic behaviour of nonlinear parabolic equations. {A}nomalous
  exponents.
\newblock In {\em Degenerate Diffusions}, volume~47 of {\em IMA series in
  Mathematics and Applications}, pages 215--227, 1993.

\bibitem{bib:wilson71a}
K.~Wilson.
\newblock Renormalization group and critical phenomena. {I. Renormalization}
  group and the {Kadanoff} scaling picture.
\newblock {\em Phys. Rev. B}, 4:3174--3183, 1971.

\bibitem{bib:wilson71b}
K.~Wilson.
\newblock Renormalization group and critical phenomena. {II. Phase}-space cell
  analysis of critical behavior.
\newblock {\em Phys. Rev. B}, 4:3184--3205, 1971.

\end{thebibliography}
\end{document}